\documentclass[a4paper]{amsart}

\usepackage{amsmath, amssymb, amsthm, graphicx, hyperref}
\theoremstyle{plain}

  \newtheorem{prop}{Proposition}
  
  \newtheorem{lem}{Lemma}
\theoremstyle{definition}

\newcommand{\cat}[1]{\mathcal{#1}}
\newcommand{\tit}{\tilde{\mathfrak{t}}}
\newcommand{\torig}{{\mathfrak{t}}}
\newcommand{\PRB}{\textit{PaRB}}
\newcommand{\FLD}{\textit{FLD}}
\newcommand{\rat}{\mathbb{Q}}

\DeclareMathOperator{\Hom}{Hom}

\begin{document}
  
\title{Formality of the chain operad of framed little disks}
\author{Pavol \v Severa}
\thanks{work supported in part by the Swiss National Science Foundation}
\address{Section de Math\'ematiques, rue du Li\`evre 2-4, 1211 Geneva, Switzerland\\
  on leave from Department of
Theoretical Physics, FMFI UK Bratislava, Slovakia}
\date{}
\maketitle

Tamarkin proved in \cite{Tam} formality of the operad of rational chains of the little disk operad. Extending his result to the framed little disk operad is a very simple exercise, and this is what we do in this note.

While writing this note I learnt about existence of a work in progress \cite{GS} containing this result and possibly much more.

All vector spaces in this note are supposed to be over $\rat$, and also all chains are with rational coefficients.

\section{Four operads}
Recall that an operad in a symmetric monoidal category $\cat{C}$ is a functor $O$ from the groupoid of finite sets $\cat{S}$ to $\cat{C}$
$$O:A\mapsto O_A$$
together with morphisms 
$$\circ_a:O_A\otimes O_B\to O_{(A-\{a\})\sqcup B}$$
for all $a\in A$; these morphisms are required to be functorial and associative in the obvious sense.

In this section we introduce the operad $BV$ of Batalin-Vilkovisky algebras and  the framed little disk operad $\FLD$. By a theorem of Getzer, the homology of $\FLD$ is the operad $BV$. The aim of this note is to show that the operad of rational chains of  
$\FLD$ is quasiisomorphic to $BV$. The quasiisomorphism is constructed via two more operads introduced in this section, namely the operad of parenthesized ribbon braids and an operad $\tit$ in the category of Lie algebras. This construction follows verbatim the paper of Tamarkin \cite{Tam}.

\subsection{Gerstenhaber and BV operad}

A Gerstenhaber algebra is a graded vector space $V$ with the following structure
\begin{enumerate}
  \item $V$ is a graded-commutative associative algebra
  \item $V[1]$ is a graded Lie algebra
  \item $[a,bc]=[a,b]c+(-1)^{(|a|+1)|b|}b[a,c]$.
\end{enumerate}

A Batalin-Vilkovisky (BV) algebra is a Gerstenhaber algebra $V$ with a linear map $\Delta:V\to V$ of degree $-1$, such that
\begin{enumerate}
  \item $\Delta^2=0$
  \item $[a,b]=(-1)^{|a|}\Delta(ab)-(-1)^{|a|}(\Delta a)b-a\Delta b$.
\end{enumerate}

Gerstenhaber operad $G$ in the category of graded vector spaces is defined in the usual way: $G_A$ is the graded vector space of all expresions built from elements of $A$ (containing each element of $A$ once in every monomial) using the operations of Gerstenhaber algebra, modulo defining relations (elements of $A$ are defined to have degree 0 in these expressions). The operad $BV$ is defined similarly.

If $A$ has $n$ elements, we have $\dim G_A=n!$ and $\dim BV_A=2^n n!$.

 \subsection{Operad of framed little disks}
  
  The framed little disk operad $\FLD$ is an operad in the category of topological spaces. $\FLD_A$ is the space of all configurations of non-overlapping disks in the unit disk, with a bijective labelling of the little disks by the elements of $A$, and with a choice of an angle for each little disk. Composition is defined by shrinking the unit disk to the size of a little disk, rotating by the corresponding angle, and glueing in. 
  
   By a theorem of Getzler \cite{Get}, the homology of $\FLD$ is isomorphic to $BV$. Namely, $\FLD_{\{a\}}$ is homotopy equivalent to $S^1$, and the generator of $H_1(S^1)$ is identified with $\Delta a$. The subspace of  $\FLD_{\{a,b\}}$ where the angles are fixed to be zero, is again homotopy equivalent to $S^1$; the generator of $H_0(S^1)$ is identified with $ab$ and the generator of $H_1(S^1)$ with $[a,b]$.

\subsection{An operad of Lie algebras $\tit$}\label{tit}
  Lie algebras form a symmetric monoidal category, with the product defined to be the direct sum.
  
  For a finite set $A$ let $\tit_A$ be the Lie algebra with generators $t_{ab}$ ($a,b\in A$, $a\neq b$) and $s_a$ ($a\in A$) and relations
  \begin{align*}
    t_{ab}&=t_{ba}\\
    [t_{ab},t_{cd}]&=0\quad\mbox{if all $a,b,c,d$ are different}\\
    [t_{ac},t_{ab}+t_{bc}]&=0\\
    [s_a,\text{anything}]&=0
  \end{align*}
  The Lie algebra $\tit_A$ is thus the direct sum of the Lie algebra $\torig_A$ generated by $t$'s with the abelian Lie algebra $\mathfrak{s}_A$ generated by $s$'s.
  
  The operadic structure on $\tit$ is defined as follows. If $A$ and $B$ are finite sets and $a\in A$ then 
  $\circ_a:\tit_A\oplus \tit_B\to \tit_{(A-\{a\})\sqcup B}$ is given by
    \begin{align*}
      \circ_a(t_{xy})=t_{xy}&\qquad\text{if }x,y\in A-\{a\}\\
      \circ_a(s_x)=s_x&\qquad\text{if }x\in A-\{a\}\\
      \circ_a(t_{xy})=t_{xy}&\qquad\text{if }x,y\in B\\
      \circ_a(s_x)=s_x&\qquad\text{if }x\in B\\
      \circ_a(t_{ax})=\sum_{y\in B} t_{xy}&\qquad\text{if }x\in A-\{a\}\\
      \circ_a(s_a)=\sum_{x\in B}s_x&\, +\sum_{\{x,y\}\subset B}t_{xy}
    \end{align*}

Out of the operad $\tit$ we can form the operad $\bigwedge\tit$ of chain complexes (computing the homology of $\tit$ with trivial coefficients). The operad $BV$ can be seen as an operad of chain complexes, with zero differential.

\begin{lem}
  There is a quasiisomorphism of operads of chain complexes $BV\to \bigwedge\tit$, given on generators of $BV$ by 
  $ab\mapsto 1$, $[a,b]\mapsto t_{ab}$ and $\Delta a \mapsto s_a$.
\end{lem}

Indeed, one readily checks that the above asignements give a morphism of operads. We need to check that it is a quasiisomorphism. Tamarkin proved it in \cite{Tam} for the morphism $G\to\bigwedge\torig$ given by $ab\mapsto 1$, $[a,b]\mapsto t_{ab}$. Since $H_\bullet(\tit_A)\cong H_\bullet(\torig_A)\otimes\bigwedge\mathfrak{s}_A$, we see that the map $BV_A\to H_\bullet(\tit_A)$ is onto and that the two spaces have the same dimension, hence it is an isomorphism.
  
  \subsection{Operad of parenthesized ribbon braids}
  Let us define an operad of groupoids $\PRB$. The objects of $\PRB_A$ are all parenthesized permutations of $A$, i.e.\ all monomials built from elements of $A$ using a noncommutative and nonassociative product (using each element once). Objects form an operad (in the category of sets) in the obvious way.

Morphisms are ribbon braids, where each ribbon connects the same elements of $A$.
  The composition $\circ_a:\PRB_A\times \PRB_B\to \PRB_{(A-\{a\})\sqcup B}$  is defined by replacing the ribbon connecting $a$'s by the $B$-ribbon braid made very thin; rotation of the $a$-ribbon is turned into rotation of the $B$-ribbon braid. 
  
  The groupoid $\PRB_A$ is obviously equivalent to the group of pure ribbon braids with $|A|$ ribbons, i.e.\ to the direct product of the pure braid group of $|A|$ strands with $\mathbb{Z}^{|A|}$.

  \section{Construction of the quasiisomorphism}
  
  \subsection{A model of the framed little disks operad}
  Out of the operad of groupoids $\PRB$ we can build the operad of its classifying spaces $|N(\PRB)|$ (the operad of geometrical realizations of its nerves). We want to see that $|N(\PRB)|$ is homotopy equivalent to $\FLD$, in the sense that there is a third topological operad $X$ and morphisms of operads
  $$\FLD\leftarrow X\to |N(\PRB)|$$
  which are homotopy equivalences.
  
  The reason is that if we choose an ordering of the set $A$, the fundamental groups of $\FLD_A$ and of $|N(\PRB)|$ become naturally isomorphic to the group $\pi_A$ of pure ribbon braids colored by elements of $A$, and that their universal covers are contractible. The spaces
  $$X_A=\textit{UC}\,(\FLD_A\times |N(\PRB_A)|\bigr)/\pi_A,$$
  where $\textit{UC}$ means universal cover,
  form the operad we looked for.
  
  Recall now that any simplicial set $K$ gives rise to a chain complex $C(K)$ (of formal linear combinations of simplices), and that the functor $C$ has a natural monoidal structure given by the Eilenberg-MacLane shuffle product. Any topological operad gives rise to the operad of its singular simplices, and  (via $C$ and its monoidal structure) to an operad of its singular chains. We now know that the operad of singular chains of $\FLD$ is quasiisomorphic (i.e.\ connected by a chain of quasiisomorphisms) with the operad of singular chains in $|N(\PRB)|$. The latter is then quasiisomorphic to $C(N(\PRB))$. In what follows we will prove that  $C(N(\PRB))$ is quasiisomorphic to its homology, which is $BV$.
  
  By abuse of notation, if $\cat{C}$ is a category, we shall denote $C(N(\cat{C}))$ simply by $C(\cat{C})$.
  
  \subsection{The morphism of Bar-Natan}
  
  An associative algebra can be seen equivalently as a linear category with one object. In this way, the enveloping algebras $U(\tit_A)$ give us an operad in linear categories, and also their completions $\hat U(\tit_A)$ (formal power series in $t$'s and $s$'s) give us such an operad.
  
  One of the results of Bar-Natan in \cite{BN} can be formulated as follows: any Drinfeld associator $\Phi$ gives rise to a morphism of operads of categories
$$\phi:\PRB\to\hat U(\tit).$$ 
On generators of $\PRB$ it is defined by

$$
\includegraphics[width=\textwidth]{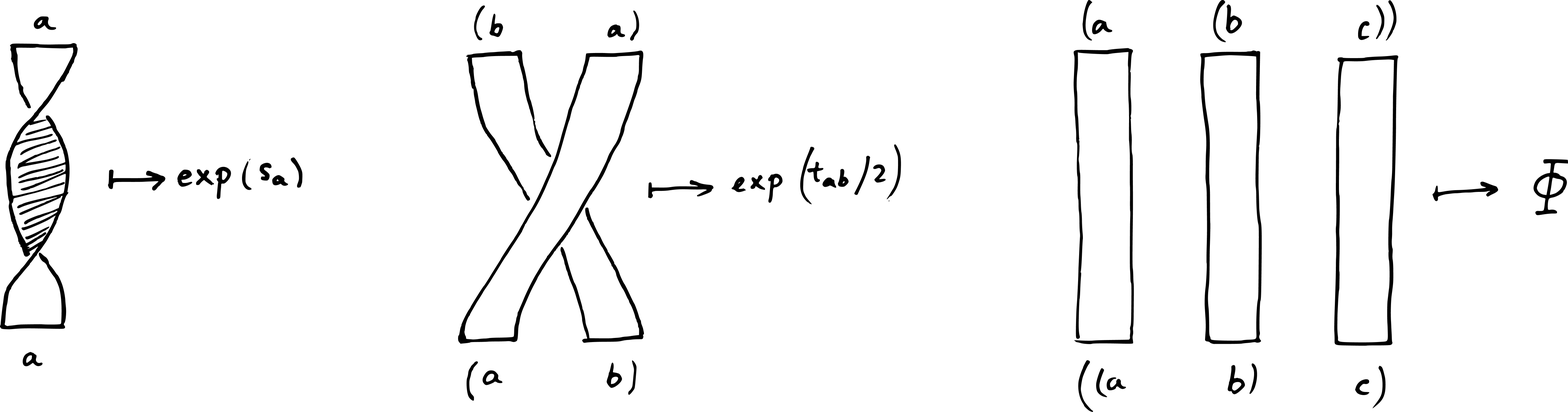}
$$

  \subsection{The linear nerve of $U(\tit)$ and its chain operad}
  
  If $K$ is a simplicial set, let $\rat K$ denote its linearization, i.e.\ the simplicial vector space with basis $K$. If $\cat{C}$ is a category, let $\rat\cat{C}$ be its linearization; it has the same objects as $\cat{C}$, and $\Hom_\cat{C}(x,y)$ is replaced by the vector space with basis $\Hom_\cat{C}(x,y)$. Let $\chi_{x,y}:\Hom_{\rat\cat{C}}(x,y)\to\rat$ be the linear map which maps all elements of $\Hom_\cat{C}(x,y)$ to $1$. It gives a functor $\chi$ from $\rat \cat{C}$ to the linear category with one object and with $\rat$ as its algebra of endomorphisms.
  
  One can easily obtain $\rat N(\cat{C})$ directly from $\rat \cat{C}$ and from the functor $\chi$. More generally, if $\cat{L}$ is a linear category and $\chi$ is a functor from $\cat{L}$ to the one object category as above, one can obviously extend this construction do define a simplicial vector space (the "linear nerve" of $\cat{L}$), denoted by $N^\rat(\cat{L})$.
  
  For any simplicial vector space $V$ we have the corresponding chain complex $C(V)$; there is a natural monoidal structure on the functor $C$ given by Eilenberg-MacLane shuffle product.
  
  Let $\chi:U(\tit_A)\to\rat$ be the counit. If we apply $N^\rat$ to the operad $U(\tit)$, we obtain an operad of simplicial vector spaces $N^\rat(U(\tit))$, and finally the operad of complexes $C(N^\rat(U(\tit)))$, denoted in what follows just by $C( U(\tit))$. The complex $C( U(\tit))$ looks as follows:
  $$C_n( U(\tit_A))\cong ( U(\tit_A))^{\otimes n}$$
  \begin{multline*}
    d(u_1\otimes\dots\otimes u_n)=\chi(u_1)u_2\otimes\dots\otimes u_n-u_1u_2\otimes\dots\otimes u_n+\dots\\+(-1)^{n-1}u_1\otimes\dots\otimes u_{n-1}\chi(u_n).
  \end{multline*}
  This is a complex computing $\operatorname{Tor}_{U(\tit_A)}(\mathbf{1},\mathbf{1})$, where $\mathbf{1}$ is the trivial $\tit_A$-representation, i.e.\ the homology of $\tit_A$ with trivial coefficients.
  
  The inclusion $\bigwedge^n\tit_A\subset(\tit_A)^{\otimes n}\subset ( U(\tit_A))^{\otimes n}$ gives a quasiisomorphism of operads $\bigwedge\tit\to C(U(\tit))$.
  The inclusion $C(U(\tit))\subset C(\hat U(\tit))$ is also a quasiisomorphism, as is easily seen from the fact that the homology $BV_A$ of $C(U(\tit_A))$ is finite-dimensional for any $A$. 
  
  \subsection{The quasiisomorphism}
  
  We already found two chains of quasiisomorphisms of operads:
  $$BV \to \textstyle\bigwedge\tit \to C(U(\tit)) \to C(\hat U(\tit))$$
  and
  $$C(\PRB) \to C^\text{sing}(|N(\PRB)|)\leftarrow C^\text{sing}(X)\to C^\text{sing}(\FLD).$$
  Now we shall connect them.

  The Bar-Natan morphism $\phi$ induces a morphism of operads $C(\PRB)\to C(\hat U(\tit))$. When we pass to homology, we get a morphism $BV\to BV$. One easily checks that it is the identity on the elements $\Delta a\in BV_{\{a\}}$, $ab\in BV_{\{a,b\}}$ and $[a,b]\in BV_{\{a,b\}}$. It is therefore the identity everywhere, and so $C(\PRB)\to C(\hat U(\tit))$ is a quasiisomorphism.
  
  We thus proved
  \begin{prop}
    The operad of rational singular chains of the framed little   disk operad is quasiisomorphic (i.e.\ connected by a chain of quasiisomorphisms of operads in the category of chain complexes) to its homology $BV$.
  \end{prop}

\end{document}